\UseRawInputEncoding
\documentclass[reqno]{amsart}

\usepackage{amsfonts}
\usepackage{bezier}
\usepackage{euscript}
\usepackage{amsfonts,amssymb,amsmath,amsbsy,amsthm,amscd}
\usepackage[english]{babel}   

\usepackage{graphics, graphicx}

\newcommand{\No}{ No }

\def\dfrac#1#2{\displaystyle{#1\over #2}}

\begin{document}
\title[On non-uniqueness in the option valuation problem
 ]{
On non-uniqueness in the option valuation problem
}


\author{Ekaterina A. Ladykova,   Olga S. Rozanova}
\address{ Mathematics and Mechanics Department, Lomonosov Moscow State University, Leninskie Gory,
Moscow, 119991, Russian Federation}
\email{rozanova@mech.math.msu.su}



\begin{abstract}
It is known \cite{Heston2006}
 that the value of a call option in the case of constant elasticity processes (CEV)  with the indicator
$\alpha$ exceeding the critical $\alpha=1$ is determined in a non-unique way.
 We show how, based on an already existing mathematical theory concerning the correctness  of boundary conditions for degenerate parabolic equations on the semi-axis $[0,\infty)$, this phenomenon can be explained. Namely,  for $1<\alpha\le \frac32$ the non-uniqueness is due to the fact that the initial data of the call option are outside the T\"acklind class, and for $\alpha> \frac32$ it is due to the absence boundary condition for $x=\infty$.
 \end{abstract}
 {
\maketitle

{\bf 1. Introduction.}
The price of an option with the exercise function $\Phi(x)$ in the case where the random process determining the dynamics of the risky asset is a process with constant elasticity (CEV), see \eqref{CEVP}, is defined as the solution of the following boundary value problem with a terminal condition [2]:
\begin{eqnarray}\label{Main}
&&U_\tau+ r x U_x + \frac12 \sigma^2 x^{2\alpha} U_{xx} - r U=0, \quad x>0,\, 0\le \tau\le T,\\
&&U(T,x)=\Phi(x).\label{MainD}
\end{eqnarray}
Here $\alpha>0$, $r\ge 0$, $\sigma>0$ are constants corresponding to the CEV process indicator, the risk-free asset trend and the volatility characteristic, respectively. In the case of the most popular ``call'' option, the execution function has the form
\begin{eqnarray}\label{call}
\Phi(x)=\max (0, x-K),\quad  K>0.
\end{eqnarray}

 The equation \eqref{Main} is best known in the case of $\alpha=1$, when the random process coincides with geometric Brownian motion, the equation in this case is named after Black and Scholes, and the formula that gives the solution to the problem \eqref{Main}, \eqref{MainD} is called the Black-Scholes formula
\cite{BS}.
It is one of the most widely used achievements of financial mathematics (e.g.,
 \cite{Hull}).
The simplicity of this formula is due to the fact that by replacing the independent variable $x=\ln y$ as a first step, the equation \eqref{Main} is reduced to the heat equation.
However, in the process of studying real markets it became clear that it is necessary to consider a wider class of random processes than geometric Brownian motion, which is associated with the interest in considering CEV processes with exponents different from one. Formulas for the option price have long been known. They were derived for all cases of the $\alpha$ exponent by reducing to the Feller equation, for which the fundamental solution is known \cite{Feller}.
For $\alpha<1$ this was done in \cite{Cox0},
and for $\alpha>1$ --- in
\cite{MacBeth}).
The formulas are quite cumbersome and contain the so-called $\chi^2$-distribution, they are included in classical textbooks, for example \cite{Hull}.
However, 25 years after the publication of the option price formula for the CEV model in \cite{Heston2006}.
it was shown that for $\alpha>1$ the option price solution is not unique, which was interpreted as the possibility of financial bubbles (a more detailed presentation can be found in \cite{Multiplicity}).
This caused a surge of interest in the problem 
among financial mathematicians who are oriented towards probabilistic methods. The non-uniqueness is explained by the fact that the CEV process for $\alpha>1$ is a strict local martingale, which leads to a violation of many familiar results of option theory 
\cite{Cox}, \cite{Heston2006}, \cite{Protter}.
A compact and neat presentation of the results of these many years of research, with derivations of formulas and corrections of inaccuracies and typos, can be found in
\cite{Hsu}, \cite{PastPresent},
however, the theory is not complete, and works on this topic continue to appear (e.g., 
\cite{Cetin}, \cite{Dias}).

We are interested in explaining the phenomenon of non-uniqueness from the point of view of the theory of linear degenerate partial differential equations.

Some aspects of this theory are not included in the circle of ideas used by specialists in theoretical financial mathematics.

In particular, we note that
traditionally, the problem \eqref{Main}, \eqref{MainD}, \eqref{call} also includes boundary conditions
\begin{equation}\label{BC_opt}
U|_{x=0}=0, \quad U|_{x=+\infty}=+\infty,
\end{equation}
natural for the option value,
while according to the theory of Fichera
\cite{Fichera}, \cite{Oleinik}
(see below) the problem requires setting boundary conditions depending on $\alpha$. Namely,
for $\alpha<\frac12$ the boundary condition is needed only at zero,
for $\alpha>\frac32$ the boundary condition is needed only at infinity,
for the remaining $\alpha$, including the Black-Scholes equation, boundary conditions are not needed, they are inherited from the final condition.

In what follows, to avoid unnecessary details, we set $r=0$ and make a change in the time variable $t=T-\tau$, after which we obtain a problem with initial conditions
\begin{eqnarray}\label{Main1}
&&U_t = \frac12 \sigma^2 x^{2\alpha} U_{xx}, \quad x>0,\, 0\le t \le T,\\
&&U(0,x)=\Phi(x).\label{MainD1}
\end{eqnarray}
In the problem of the price of a call option, the initial condition obtained in this way has a linear growth at infinity, but for the study of the phenomenon of non-uniqueness we will not limit ourselves to such special initial data.

Thus, in light of the results on the rules for setting boundary conditions in accordance with Fickera's theory, the following problems arise:

1. On the uniqueness class for $\alpha<\frac12$ with growth constraint at infinity (Tacklind class);

2. On the uniqueness class for $\alpha>\frac32$ with growth constraint at zero;

3. On the uniqueness class for $\frac12\le\alpha\le \frac32$ with growth constraint at zero and at infinity.
In particular, we will show that the non-uniqueness of the solution in the option pricing problem for $\alpha>1$ occurs for different reasons: for $1<\alpha\le\frac32$ the initial data of linear growth do not belong to the T\"acklind class, and for $\alpha>\frac32$ the problem lacks a boundary condition at infinity. Apparently, such a separating role of the exponent $\alpha=\frac32$ in this context was not known, although in other problems the value of the exponent $\frac32$ also arises. For example, in \cite{Meyer}
(section 1.2.2) the value of this exponent is indicated as separating when setting boundary conditions at infinity, and in \cite{Heston2006}
it is noted that the CEV process has finite variance for $\alpha>1$ only for $\alpha>\frac32$.

Note that the question of the growth of the solution at infinity, due to the replacement of the independent variable $x=\frac{1}{z}$, turns out to be connected with the question of the growth of the solution at zero (for other exponents $\alpha$), see Sec. 4.

Also note that the equation \eqref{Main1} by replacing $x=\left(\frac{y}{1-\beta}\right)^{1-\beta}$, $\alpha = \frac{\beta}{1-\beta}$ can be reduced to an equation containing the Bessel operator $B$ of the form
\begin{equation*}
U_t=B_y U, \quad B_y = y^{-\beta} \left(\dfrac{\partial}{\partial y} y^\beta \dfrac{\partial}{\partial y}\right) = \dfrac{\partial^2}{\partial y^2} + \frac{\beta}{y} \dfrac{\partial}{\partial y}.
\end{equation*}
Such an equation is called the ``$B$-equation'' \cite{Kiprianov}
or the ``heat  equation on a fractal'' \cite{Nahushev}.
\bigskip

{\bf 2. Fundamental solution.}
The equation \eqref{Main1} (as well as \eqref{Main}) by means of the transformation
$$y=\frac{1}{\sigma^2 (1-\alpha)^2} x^{2(1-\alpha)}$$ is reduced to the equation
\begin{equation}\label{FellerEq}
u_t\,=\,(a\,y\,u)_{yy}\,-\,((b\,y\,+\,c)\,u)_y,\quad y>0,
\end{equation}
where $a, b, c$ are constants, $a\,>\,0$,
\begin{align*}\label{ds16}
a\,=\,{2},\quad
b\,=\,0,\quad 
c\,=\,\frac{3-2\,\alpha}{1-\alpha}.
\end{align*}

Feller \cite{Feller}
constructed a theory of equations of the form \eqref{FellerEq}. In particular, its
fundamental solution has been found, that is, the solution \eqref{FellerEq} with the initial conditions
\begin{equation*}
u|_{t=0}= \delta(y-\xi).
\end{equation*}

We write it in the form given in \cite{PastPresent}
for $b=0$, where the typo in the original work has been corrected:
\begin{equation}\label{E}
{\mathcal E}_\alpha(t,\,y;\,\xi)\,=\frac{1}{2t}\,\left(\frac{y}{\xi}\right)^{\nu/2}\,\exp\left(-\frac{y+\xi}{2t}\right)
\,I_{|\nu|}\left(\frac{\sqrt{y\xi}}{t}\right),
\end{equation}
here
$I_k(x)\,=\,\sum\limits_{r\,=\,0}^{\infty}\frac{(x/2)^{2r\,+k}}{r!\,\Gamma(r\,+\,1\,+\,k)}$
-- modified
Bessel function of the first kind,
$\nu= 1-\frac{c}{a}=-\frac{1}{2(1-\alpha)}$. This expression is valid for all values ??of $\alpha\ne 1$.
Note that the fundamental solution in the case $b=0$ was obtained by S.~Kepinski in 1905
long before Feller's work.

Since $I_k(z)=O(z^k)$, $z\to 0$, then
\begin{itemize}
\item function ${\mathcal E}_\alpha$ is finite at $y=0$, namely:

- at
$0<\alpha<1$ the value of ${\mathcal E}_\alpha(t,\,0;\,\xi)$ is nonzero,

- at $\alpha >1$ the value of ${\mathcal E}_\alpha(t,\,0;\,\xi)=0$ and ${\mathcal E}_\alpha(t,\,y;\,\xi)=O(y^\nu)$, $y\to 0$.

\item For $y\to \infty$, it follows from the asymptotics of the Bessel function that
$${\mathcal E}_\alpha(t,\,y;\,\xi)\sim c_1 \,e^{-c_2 y} y^{\frac{2\nu -1}{4}} \cos\left(c_3 \sqrt{y}-\left(|\nu|+\frac12\right) \frac{\pi}{2}\right),$$
$c_i$ are positive constants that depend on $t$ and $\xi$.
\end{itemize}For $0<\alpha<1$ ($\nu<0$), the integral over the positive semi-axis of the solution given by \eqref{E} is not preserved (decreases in time)
\cite{PastPresent}:
\begin{equation*}
\int\limits_0^\infty {\mathcal E}_\alpha(t,\,y;\,\xi) \, dy =\Gamma\left(-\nu, \frac{\xi}{2t}\right)<1, \quad \Gamma(m,z)=\frac{1}{\Gamma(m)}
\int\limits_0^z \eta^{m-1} e^{-\eta} d\eta,
\end{equation*}
therefore, in order to preserve the integral, it is necessary to add a singular component. Thus, the fundamental solution ${\mathcal E}(t,\,y;\,\xi)$ has the form
\begin{equation*}
{\mathcal E}(t,\,y;\,\xi)={\mathcal E}_\alpha(t,\,y;\,\xi)+ \left(1-\Gamma\left(-\nu, \frac{\xi}{2t}\right)\right) \delta(\xi).
\end{equation*}
For $\alpha >1$, from the asymptotic properties of ${\mathcal E}_\alpha(t,\,y;\,\xi)$ at zero and at infinity it follows that the total flux
$\left[(a y{\mathcal E})_y-c{\mathcal E}\right]\Big|_0^\infty=0 $, therefore the integral over the positive semi-axis of the solution is preserved.

In the limit as $t\to \infty$ for each fixed positive $y$ and $\xi$ the function ${\mathcal E}_\alpha(t,\,y;\,\xi)$ tends to zero. This follows, for example, from the non-existence of a stationary solution of the equation \eqref{FellerEq} bounded for $y>0$ (the general solution is $C_1+C_2 y^{-\nu}$).

The expression for the fundamental solution $\bar {\mathcal E}_\alpha(t, x; \eta)$ for \eqref{Main1} follows from \eqref{E} after the replacement
\begin{equation*}
y=\frac{4 \nu^2}{\sigma^2} x^{-\frac{1}{\nu}}, \qquad \xi=\frac{4 \nu^2}{\sigma^2} \eta^{-\frac{1}{\nu}}.
\end{equation*}

The properties of \eqref{E} imply the following properties of $\bar {\mathcal E}_\alpha$.

1. For $0<\alpha<1$:
\begin{itemize}
\item $\bar {\mathcal E}_\alpha$ is finite for $x=0$ and $\bar {\mathcal E}_\alpha= O(\exp(-C x^{2(1-\alpha)})$, $x\to\infty$, $C>0$;
\item $\bar {\mathcal E}_\alpha=O(\eta)$, $\eta\to 0$ and $\bar {\mathcal E}_\alpha= O(\exp(-C \eta^{2(1-\alpha)})$, $\eta\to\infty$, $C>0$.
\end{itemize}

2. For $\alpha>1$ the origin and infinity change in places:
\begin{itemize}
\item $\lim\limits_{x\to 0} \bar {\mathcal E}_\alpha(t,\,x;\,\eta)=0$, and $\bar {\mathcal E}_\alpha(t,\,x;\,\eta)=O(\frac{1}{x}),$ $x\to\infty$ (from here, in particular, follows the conditional convergence of the integral of $\bar {\mathcal E}_\alpha(t,\,x;\,\eta)$ along the axis $x>0$ and the change of sign of $\bar {\mathcal E}_\alpha(t,\,x;\,\eta)$);

\item $\lim\limits_{\eta \to 0} \bar {\mathcal E}_\alpha(t,\,x;\,\eta)=0$, and $\bar {\mathcal E}_\alpha(t,\,x;\,\eta)=O(1),$ $\eta\to\infty$. It follows that in order for the solution of the problem \eqref{Main1}, \eqref{MainD} to be found as a convolution of the fundamental solution with the initial data, it is necessary that the condition $\Phi(x)=o(x)$, $x\to \infty$ be satisfied. This prevents the possibility of calculating the value of the call option with the final data \eqref{call}.
\end{itemize}

Note that the expression for the transition density for \eqref{Main} obtained in \cite{Cox0}, \cite{MacBeth} follows from the analogue of \eqref{E} for the more general case with $b\ne 0$.

The solution of the Cauchy problem \eqref{Main1}, \eqref{MainD1} is found as
\begin{equation}\label{SolE}
U(t,x)=\int\limits_0^\infty \Phi(\eta) \bar {\mathcal E}(t,\,x;\,\eta)\,d\eta
\end{equation}
provided that this integral converges.
\bigskip

{\bf 3. Exact solutions.}
Note that some classes of solutions of the equation \eqref{FellerEq} are known. We present them for the case $b=0$. This is the stationary solution
\begin{equation*}
u=C_1+C_2 \,y^{-{\nu}}, \quad C_1, C_2 = {\rm const},
\end{equation*}
and the solution with separable variables
\begin{eqnarray}\label{sep}
u&=& u_1(t) u_2(y) \\
u_1 &=& K e^{-kt}, \quad k>0,\quad u_2=C_1\, y^{-\frac{\nu}{2}}\,J(-\nu, \sqrt{2c y})+C_2 \,y^{-\frac{\nu}{2}}\,Y(-\nu, \sqrt{2c y}),\quad C_1, C_2 = {\rm const},\nonumber
\end{eqnarray}
where $J$ and $Y$ are Bessel functions of the 1st and 2nd kind. We see that for every fixed $x>0$ the solution \eqref{sep} tends to zero as $t\to \infty$.

In addition, in \cite{Feller}
another exact solution is given, which we will write out in the case $b=0$:
\begin{equation}\label{sol}
u(t,x)=t^{\frac{c}{a}-2}\, \exp\left(-\frac{y}{a t} \right), \quad \frac{c}{a}-2=-\nu-1,
\end{equation}
satisfying the zero initial conditions for $x>0$. In \cite{Feller}
it is noted that if 
$c<a$ ($\alpha<1$), then the integral $\int\limits_0^y u(t,z) dz$ diverges as $t\to 0$, which indicates the presence of a non-integrable singularity at zero.

The solutions of the equation \eqref{Main1} are obtained from here by replacing $y=\frac{4 \nu^2}{\sigma^2} x^{-\frac{1}{\nu}}$. In particular, we see that the stationary solution has the form $U(t,x)=C_1+C_2 \,x $. As for the solution \eqref{sep}, its $x$-dependent component for $\alpha<1$
is oscillating with an amplitude increasing at infinity, and for $\alpha>1$ at infinity the solution is equivalent to $C_1+C_2 \,x $ (this follows, for example, from \cite{Bellmann},
Theorem 5, Chapter 6).

\bigskip
According to \cite{Oleinik}
the solutions belong to the class of bounded measurable functions and are understood in the weak sense (in the sense of an integral identity), the smoothness of such solutions is very delicate \cite{Castro}.
2. Feller in \cite{Feller}
classified the boundary $y=0$ for the equation \eqref{FellerEq} from the point of view of the motion of a particle subject to a stochastic process
\begin{equation}\label{FP}
dY=(bY+c) dt+\sqrt{2 a Y} dW,
\end{equation}
where $dW$ is the differential of the standard Brownian motion. The density of the process \eqref{FP} is described by the equation \eqref{FellerEq}.
The boundary can be attainable (reflecting or absorbing) or unattainable. For $c\le 0$ it is attainable and absorbing, for $0< c< a$ the boundary is attainable, it can be reflecting or absorbing, for $ c> a$ the boundary is unattainable. Thus, boundary conditions are required only for the case $0< c< a$.

The CEV process is described by the stochastic equation
\begin{equation}\label{CEVP}
 dF=\sigma F^\alpha dW,
\end{equation}
which by transforming $Y=\frac{1}{\sigma^2 (1-\alpha)^2} F^{2(1-\alpha)}$ can be converted to \eqref{FP} when
\begin{align*}\label{ds16}
 a\,=\,{2},\quad
b\,=\,0,\quad 
c\,=\,\frac{1-2\,\alpha}{1-\alpha}.
\end{align*}
Thus, for the process \eqref{CEVP} the boundary $f=0$ is attainable (reflecting or absorbing) for $\alpha< \frac12$,
it is attainable and absorbing for $\frac12\le \alpha <1$, and the boundary is unattainable for $\alpha >1$.
This classification corresponds to the rules for setting boundary conditions from the point of view of Fichera's theory.

\bigskip

{\bf 5. On Tikhonov-T\"aklind classes.}
Questions of uniqueness of the solution of the Cauchy problem for linear parabolic equations have been studied since Holmgren's 1924 work on the heat equation. However, the fundamental work should be considered Tikhonov's 1935 work
\cite{Tikhonov},
in which
not only was the uniqueness class indicated
\begin{equation*}\label{heat}
|u(х, t)|\le B \exp {\beta |x|^2 }, \quad \beta>0, \, x \in {\mathbb R}, \,t\le 0,
\end{equation*}
but its optimality was demonstrated. In \cite{Taklind}
Taklind refined Tikhonov's result for the heat equation, showing that the uniqueness of the solution of the Cauchy problem
occurs in the class
\begin{equation}\label{heat}
|u(х, t)|\le B \exp{|x| h(x) }, \quad x \in {\mathbb R}, \,t\le 0,
\end{equation}
where $h(x)$ is a non-decreasing non-negative function such that $\int\limits_1^\infty \frac{1}{h(x)} dx=\infty$.
The literature in which these results are extended to equations with variable coefficients is quite extensive (see the classic review \cite{Ilyin}
).

\medskip

{\bf 5.1 $\alpha<1$, $x\to \infty$, unattainable boundary.}
In this case, the boundary condition at infinity is not required. The results of \cite{Kamynin_Khimchenko}, \cite{Kamynin} are applicable to study the uniqueness classes for functions increasing at infinity in the problem \eqref{Main1}, \eqref{MainD}.
Namely, from Theorem 1 \cite{Kamynin_Khimchenko}
it follows that if $g(s)$ is such that
$s^{2\alpha}\le g^2(s)$, $s>s_0>0$, and $\int\limits_{s_0}^\infty \frac{1}{g(s)} ds=\infty$, then the solution is unique in the class of functions
\begin{equation}\label{T1}
|u(t, x)|\le B \exp{(G(x) h(G(x)) )},\quad G(x)=\int\limits_{s_0}^x \frac{1}{g(s)} ds \quad x \in {\mathbb R}, \,t\le 0,
\end{equation}
$h(s)$ is the T\"acklind function.
The rigorous formulation of this theorem (which is quite cumbersome), applicable in a much more general case than ours, is not given here, referring to the original paper.

It is easy to see that $g(s)=s^\alpha$, and only in the case $\alpha\le 1$ does the integral diverge. In this case $G(x)= C\,x^{1-\alpha}$, $C>0$ if $\alpha<1$ and $G(x)= C\,|\ln x|$ if $\alpha=1$ (which corresponds to the original Black-Scholes equation).

The estimate \eqref{T1} is optimal. It is easy to see that this result overlaps the much later result on the uniqueness of the solution of the Cauchy problem for the equation \eqref{Main1} in the class of functions of sublinear
growth \cite{Convexity},
Theorem 4.3, see below.

{\bf 5.2. $\alpha>1$, $x\to 0$, unattainable boundary.}
In this case, the boundary condition at zero is not required.
Since by replacing $x\to \frac{1}{x}$ we can transfer the infinity point to zero, then from the results of the previous section for $x\to 0$ we obtain the following growth constraints on the solution at zero:
\begin{equation*}\label{T2}
|u(x, t)|\le B \exp{( x^{1-\alpha} h(x^{1-\alpha})) },\quad x\in {\mathbb R}, \,t\le 0.
\end{equation*}

{\bf 5.3 $\alpha>1$, $x\to \infty$.}
In this case, for $\alpha>\frac32$ the boundary condition is required, but for $1<\alpha\le\frac32$ it is not. Important results concerning all values ??of the exponent $\alpha$ are contained in \cite{Convexity}.
Theorem 3.2 shows that
\begin{quote}\it
for all $\alpha\ge 0$, if the initial condition \eqref{MainD1} is continuous and has at most linear growth at infinity, then the solution \eqref{MainD1}, found as a convolution with the fundamental solution, has at most linear growth.
\end{quote}
Note that in the case $\alpha>1$ the integral \eqref{SolE} expressing this convolution is not necessarily convergent due to the absence of decay of the fundamental solution at infinity $x\to\infty$, however, as shown in the proof of the theorem, a finite limit of integrals over a finite interval $(0,M)$, $M\to \infty$, does exist. Moreover, as follows from Theorem 4.1,
\begin{quote}\it
the solution found as \eqref{SolE}, with initial data having linear growth at infinity (like \eqref{call}), at any time $t>0$ is a function growing slower than any positive power of $x^\lambda$, $\lambda>0$, and for $\alpha>\frac32$ it is bounded.
\end{quote}

However, such a solution is not unique. Indeed, the function $U_{st}=x$ is a time-independent solution of the equation \eqref{Main1} for all $\alpha$, but if for $\alpha\le 1$ the convolution of the fundamental solution with the initial data $\Phi(x)=x$ yields this solution (we are in the uniqueness class), then for $\alpha>1$, as follows from \cite{Convexity},
Theorem 4.1, we obtain another solution (respectively, infinitely many solutions).

An explicit example of such another solution can be constructed for $\alpha=2$:
\begin{equation*}
U_b(t,x)= x \left(1-2 {\bf\Phi}\left(-\frac{1}{\sigma x \sqrt{t}}\right)\right),
\end{equation*}
where $ {\bf\Phi}$ is the normal distribution function.
It is easy to see that this solution is bounded for $t>0$ (in \cite{Cox}
the example is attributed to Johnson, Helms, see also \cite{Multiplicity}
for a derivation of this result).

In \cite{Convexity},
Theorem 4.3, it is also proved that
\begin{quote}\it
in the class of functions $U(t,x)$ of sublinear growth in $x$, that is,
for functions such that
\begin{equation*}\label{sublin}
|{U(t,x}|\le f(x), \, t\ge 0,\quad
\lim\limits_{x\to\infty} \frac{f(x)}{x}=0,
\end{equation*}
the solution to the Cauchy problem \eqref{Main1}, \eqref{MainD1} is unique.
\end{quote}

It is important to note that in this theorem the uniform boundedness of the class of solutions with respect to $t$ is important. Indeed, it is easy to verify that the function
\begin{equation}\label{solx}
U_\nu (t,x)=t^{-\nu-1}\, \exp\left(-\frac{2 x^{-\frac{1}{\nu}}}{\sigma^2 \nu^2 t} \right)
\end{equation}
is a solution of the equation \eqref{Main1} for all values ??of $\nu\ne 0$ ($\alpha\ne 1$), the function is bounded for each $t$ and
$\lim\limits_{t\to 0} u(t,x)=0$, $x\ge 0$. This example is constructed based on the solution of \eqref{sol}. A comparison with the identically zero solution would seem to show that even in the class of bounded functions the solution of the Cauchy problem is not unique. However, $\sup\limits_{x\in \mathbb R} U_\nu\to\infty$, $t\to 0$, although for every fixed $x\in \mathbb R $ we have $ U_\nu(t,x)\to 0$, $t\to 0$.

{\bf 5.4. $\alpha<1$, $x\to 0$.}
In this case, for $\alpha<\frac12$ the boundary condition is required, but for $1>\alpha\ge\frac12$ it is not.

It follows from the results of the previous section that the uniqueness of the solution of the Cauchy problem for $x\to 0$ is ensured by the uniform in $t$ condition
\begin{equation*}\label{superlin}
|{U(t,x}|\le f(x), \, t\ge 0,\quad
\lim\limits_{x\to 0} x{f\left(\frac{1}{x}\right)}=0.
\end{equation*}
Note that the solution may have an integrable singularity at zero (for example, of the form $x^{\epsilon-1}, \, \epsilon>0$).
A solution \eqref{solx} with the property $\lim\limits_{t\to 0} u(t,x)=0$, $x\ge 0$ and bounded for each $t$, but not uniformly,
$\sup\limits_{x\in \mathbb R} U_\nu\to\infty$, $t\to 0$, has a non-integrable singularity at zero (a singular component, as noted in p.2).
\bigskip

{\bf 6. Example of the original Black-Scholes equation.}
For the Black-Scholes equation $\alpha=1$. In equation \eqref{Main1} we make a change of variable $y=\ln x-\frac12 \sigma^2 t$, and for the function $U(t, x(y,t))=\bar U (t,y) $ we obtain the heat equation defined on the entire axis:
\begin{equation*}
\bar U_t=\frac12 \sigma^2 \bar U_{yy}, \quad y \in \mathbb R.
\end{equation*}
We use classical results on Tacklind classes for the heat equation and obtain from \eqref{heat} the following uniqueness class:
\begin{equation*}\label{BS}
|u(х, t)|\le B \exp{( |\ln (x)| h(\ln (x)) )}, \quad x \in {\mathbb R}, \,t\le 0,
\end{equation*}
where $h$ is the T\"acklind function. We see that the same result follows from \eqref{T1}. We see that the initial data corresponding to the value of the ``call'' option are certainly in this class. As for the boundary conditions \eqref{BC_opt}
at zero and infinity, they should not be prescribed in advance, but follow exclusively from the initial data. The solution to such a problem is unique.

\bigskip

{\bf 7. Conclusion.}
From the above analysis it follows that the non-uniqueness of the solution in the problem \eqref{Main1}, \eqref{call} for $\alpha>1$ occurs for the following reasons.

1. For $1<\alpha\le \frac32$ for $x=+\infty$ no condition should be imposed, but the solution is unique only in the class of functions that have sublinear growth at infinity, where the initial data \eqref{call} do not fall.

2. For $\alpha> \frac32$ an additional condition on the behavior of the solution for $x=+\infty$ is required, otherwise the problem is underdetermined.

Note that in \cite{Heston2006}
methods for isolating a unique solution are proposed, for example, through the parity of the ``call'' and ``put'' options.

Taking into account the results described above, we see that uniqueness can be achieved, for example, by requiring that the solution be found using the formula \eqref{SolE}. Note that this automatically imposes the requirement that the solution be bounded for $\alpha> \frac32$.

We emphasize that we do not claim that most of the presented results are new (only some examples and interpretations are new). We only bring them together, analyze them, and offer a new perspective on the problem.

\end{document}